\documentclass[11pt]{article}

\usepackage{cite}
\usepackage{amsmath}
\usepackage{amsfonts}

\def\corresponds{{\lower.2ex\hbox{=}}{\rm\kern-.75em^\triangle}}
\def\succsim{\succ\kern-.9em_\sim\kern.3em}
\def\precsim{\prec\kern-1em_\sim\kern.3em}
\def\slantfrac#1#2{\kern1em^{#1}\kern-.3em/\kern-.1em_{#2}}
\def\lfrac#1#2{{}^{#1\!}\kern-.0em/_{#2}}
\def\buildrel#1\under#2{\mathrel{\mathop{\kern0pt #2}\limits_{#1}}}

\newcommand{\Ent}[1]{[\mkern - 2.5 mu [#1] \mkern - 2.5 mu ]}

\begin{document}

\bibliographystyle{myprsty}

\vspace{2cm}

\noindent
\null \hfill \null \\[3ex]
\vspace{0.5cm}

\begin{center}
\begin{tabular}{c}
\hline
\rule[-5mm]{0mm}{15mm} 
{\Large \sf Application of the
Combined Nonlinear--Condensation Transformation} \\
\rule[-5mm]{0mm}{10mm}
{\Large \sf to Problems in Statistical Analysis and Theoretical Physics}\\
\hline
\end{tabular}
\end{center}
\vspace{0.5cm}
\begin{center}
Sergej~V.~Aksenov$^{a)}$, Michael~A.~Savageau$^{a)}$,
\\[1ex]
Ulrich~D.~Jentschura$^{b),c)}$, Jens~Becher$^{b)}$, Gerhard~Soff$^{b)}$,
\\[1ex]
and Peter J.~Mohr$^{d)}$ 
\end{center}
\vspace{0.05cm}
\begin{center}
$^{a)}${\it Department of Microbiology and Immunology,\\
University of Michigan, Ann Arbor, MI 48109, USA}\\[1ex]
$^{b)}${\it Institut f\"{u}r Theoretische Physik,}\\
{\it Technische Universit\"{a}t Dresden,
01062 Dresden, Germany}\\[1ex]
$^{c)}${\it Theoretische Quantendynamik, Fakult\"{a}t 
Physik der Universit\"{a}t Freiburg,}\\
{\it Hermann--Herder--Stra\ss{}e 3, 79104 Freiburg, Germany}\\[1ex]
$^{d)}${\it National Institute of Standards and Technology,\\
Mail Stop 8401, Gaithersburg, MD 20899, USA}
\end{center}

\vspace{0.6cm}

\begin{center}
\begin{minipage}{11.5cm}
{\underline{Abstract}}
We discuss several applications of the recently proposed combined
nonlinear--condensation transformation (CNCT) for the evaluation of slowly
convergent, nonalternating series. These include 
certain statistical distributions which are of importance in
linguistics, statistical-mechanics theory, and biophysics
(statistical analysis of DNA sequences).
We also discuss applications of the transformation in
experimental mathematics, and we briefly expand on further
applications in theoretical physics.  Finally, we discuss
a related Mathe\-matica program for the computation of 
Lerch's transcendent. 
\end{minipage}
\end{center}

\vspace{0.6cm}

numerical approximation and analysis; DNA, RNA; \\
Folding and sequence analysis.
\vfill
\begin{center}
\begin{minipage}{16cm}
\begin{center}
\hrule
{\bf \scriptsize
electronic mail: aksenov@umich.edu and jentschura@physik.uni-freiburg.de.}
\end{center}
\end{minipage}
\end{center}



\newpage

%
%
\section{Introduction}

The combined nonlinear--condensation
transformation~\cite{JeMoSoWe1999} (CNCT) is an algorithm that
transforms a sequence of partial sums of a slowly convergent
nonalternating (monotone) series (i.e., a series whose terms all
have the same sign) into a sequence of transforms (to be referred to
as CNC transforms) often with better numerical properties. A
number of applications (e.g.~\cite{JeMoSo1999,JeGiVaLaWe2001}) have
recently been described; an acceleration of the convergence by
several orders of magnitude is observed in many cases. The CNCT
addresses, at least in part, the difficulties usually associated
with the acceleration of the convergence of monotone series. The
application of convergence acceleration methods to series of this
type typically leads to numerical instabilities and cancellations in
higher transformation orders (for a discussion see~\cite{We1989}).
In many cases, it has been observed that the CNCT is a remarkably
stable numerical process, mainly because the nonalternating input
series is transformed into an alternating series {\em before}
the actual convergence acceleration method is applied. The CNCT is a
process involving two {\em \'{e}tappes}; first step: transformation
{\em nonalternating}$\to${\em alternating series}, second step:
convergence acceleration via the delta
transformation (this transformation was introduced
in~\cite{Si1981}, and its usefulness as a powerful generalized resummation and
convergence acceleration method was demonstrated in~\cite{We1989}).

The purpose of this paper is threefold: 
(i) to discuss applications of the CNCT in such diverse disciplines 
as statistics, theoretical physics, and experimental mathematics, 
(ii) to propose an implementation of the CNCT for the
calculation of Lerch's transcendent in Mathematica~\cite{Wo1988,disclaimer},
and (iii) to provide numerical evidence for the convergence of the
CNCT in higher transformation order using computer arithmetic with
enhanced accuracy~\cite{Ba1990tech,Ba1993,Ba1994tech}.

Numerical evidence gained from high-precision calculations might be useful
because, currently, no direct, general proof exists for the convergence of
the delta transforms~\cite{Si1981,We1989}, which form part of the second
step of the CNCT. Of course, it should be noted that the convergence of
the delta transforms has been proven for a set of model
problems~\cite{Si1981}, and the exactness of the CNC transformation has
been shown in~\cite{JeMoSoWe1999} for a number of model series, but this
information is of limited use for the investigation of realistic
applications. This fact has to be contrasted with the observed rapid
convergence of the transforms in many applications of practical importance
(e.g.~\cite{JeMoSo1999,JeGiVaLaWe2001}). In particular, the recent
evaluation~\cite{JeMoSo1999} of the bound-electron self-energy has lead to 
results with
roughly 19 significant figures for atomic hydrogen (the CNCT has been used
in this calculation for the acceleration of slowly convergent partial wave
expansions). For atomic hydrogen, the renormalization process leads to a
numerical loss of about 9 significant figures, because the physically
significant part of the bound-electron self-energy is obtained after
subtracting the self-energy of the free electron, the latter is a part of
the electron mass. The succesful verification of the consistency of the
remaining 10 significant figures left after the renormalization with
higher-order analytic results~\cite{Pa1993} 
provides a sensitive test of the
numerical results and, consequently, of the convergence properties and the
reliability of the algorithm described in~\cite{JeMoSoWe1999}
and in the 
current article.
To complement this evidence, we investigate here the CNCT with
extended-precision arithmetic and gauge the rate of convergence of the CNC
transforms in higher order.

This paper is organized as follows: in Sec.~\ref{sec_cnct}, we
recall the motivation and the formulas used in the construction of
the CNCT. In
Sec.~\ref{sec_stat}, we discuss the application of the CNCT to the
evaluation of a class of important special functions,
which may be used e.g.~for the statistical analysis of DNA sequences. 
In this context we also discuss an explicit
implementation of the CNCT in Mathematica
which may serve as a template for
applications in other areas and programming languages.
These investigations are supplemented by the
high-precision calculations in Sec.~\ref{sec_math}, where a certain
application is discussed which could probably best be associated to
the field of ``experimental mathematics''. 
In Sec.~\ref{sec_other},
we discuss further applications of the CNCT, mainly related
to the evaluation of other special functions which occur naturally
in theoretical physics. 
We provide concluding remarks in Sec.~\ref{conclu}.  
Appendix~\ref{app} contains a detailed description of 
Mathematica code for the 
calculation of Lerch's transcendent.

%
%
\section{The Combined Nonlinear--Condensation Transformation}
\label{sec_cnct}

%
%
\subsection{CNCT: Motivation}
\label{cnct_motivation}

We briefly recall the notion of logarithmic convergence and the
difficulties associated with the acceleration of the convergence of
nonalternating series. Let a sequence $\{ s_n\}_{n=0}^{\infty}$
fulfil the asymptotic condition
\begin{equation}
\label{LimRatSeq}
\lim_{n \to \infty} \frac {s_{n+1} - s} {s_n - s} \; = \; \rho \,,
\end{equation}
where $s = s_{\infty}$ is the limit of the sequence as $n \to
\infty$. If $\rho > 1$, then the sequence $\{
s_n\}_{n=0}^{\infty}$ is divergent. For $\rho = 1$, the
sequence may either be convergent or divergent. A convergent
sequence with $\vert\rho\vert = 1$ is called {\em logarithmically
convergent} (if $\rho < 1$, the series is called {\em
linearly convergent}). Let us further assume that the elements of
the sequence $\{ s_n \}_{n=0}^{\infty}$ in Eq.\ (\ref{LimRatSeq})
represent partial sums 
\begin{equation}
\label{defsn}
s_n = \sum_{k=0}^{n} a(k)
\end{equation}
of an infinite series. Here, we will almost exclusively investigate
slowly convergent nonalternating sequences $\{ s_n\}_{n=0}^{\infty}$
whose elements are all real and positive (for these sequences, $0 <
\rho \leq 1$). In the case of slow convergence, $\rho$ is either
very close or equal to unity.

As observed by many authors~(e.g.~\cite{We1989}), the acceleration
of the convergence of nonalternating sequences is a potentially
unstable numerical process. The reason is the following: A sequence
transformation can only accelerate convergence if it succeeds in
extracting additional information about the index-dependence of the
remainders or ``truncation errors''
\begin{equation}
\label{defrn}
r_n = s_n - s 
\end{equation} 
from a necessarily finite set of partial sums $\{ s_n \}_{n=0}^k$ of
the input series. Normally, this is done by forming arithmetic
expressions involving higher weighted differences of the $s_n$. The
calculation of higher weighted differences is a potentially unstable
process which can easily lead to a serious loss of numerical
significance due to cancellation if the input data all have the same
sign.

The main notion of convergence acceleration is to
extract information ``hidden in trailing digits''
from a necessarily finite number of sequence elements,
in order to convert a sequence $\{s_n\}_{n=0}^{\infty}$
into a new sequence $\{s'_n\}_{n=0}^{\infty}$ with hopefully
better numerical properties. 

Let us assume
that the input sequence $\{ s_n \}_{n=0}^\infty$
is logarithmically or linearly convergent (with $\rho \leq 1$) and 
that a sequence transformation can be constructed so that
the transformed sequence
$\{ s'_n \}_{n=0}^\infty$ is linearly convergent with 
$\rho' < \rho$
[see Eq.~(\ref{LimRatSeq}) for the definition of $\rho$].  
Convergence of the new sequence can be said to be accelerated if $\rho' < \rho$.  
An equivalent definition of convergence acceleration can 
be given by calculating the limit of the following sequence
of the ratios $\chi_n$,
\begin{equation}
\chi_n \equiv \frac{s'_n - s}{s_n - s}\,,
\end{equation}
and to {\em define} convergence to be {\em accelerated} if 
\begin{equation}
\label{ConvAccelDef}
\lim_{n \to \infty} \chi_n = 0\,.
\end{equation}
Indeed, if $\rho' < \rho$, then
\begin{equation}
\label{ConvAccelDef2}
\lim_{n \to \infty} \frac{\chi_{n+1}}{\chi_n} \equiv 
\lim_{n \to \infty} \frac{s'_{n+1} - s}{s'_n - s} \,
\frac{s_n - s}{s_{n+1} - s} = \frac{\rho'}{\rho} < 1\,.
\end{equation}
So, if $\rho' < \rho$, then $\chi_n$ asymptotically behaves as a geometric 
progression within its circle of convergence,
and thus $\lim_{n \to \infty} \chi_n = 0$.

%
%
\subsection{CNCT: Formulas}
\label{cnct_formulas}

Recently, the combined nonlinear--condensation
transformation~\cite{JeMoSoWe1999} (CNCT) has been proposed as a
computational tool for the accelerated numerical evaluation of
slowly convergent, nonalternating series. The idea is to divide the
acceleration process in two steps. The first step, which is a
reordering process, consists in a rearrangement of the terms of the
input series into an alternating series via a Van Wijngaarden
transformation~\cite{vW1965}. The output of the
first step is an alternating series whose terms do {\em not}
decay more rapidly in magnitude than those of
the original input series (see the discussion in Sec.~3
of~\cite{JeMoSoWe1999}).
It could appear that nothing substantial has been
achieved in the first step of the CNCT. The Van Wijngaarden step of the CNCT
merely represents a ``computational investment'' with the intention
of transforming the nonalternating input series into a form which
is more amenable to the acceleration of convergence. The second
step, which represents a convergence acceleration process, consists
in the application of a powerful nonlinear sequence transformation
for the acceleration of the convergence of the {\em alternating}
series which resulted from the first step of the CNCT.

Following Van Wijngaarden \cite{vW1965}, we transform the
nonalternating input series 
\begin{equation}
\sum_{k=0}^{\infty} a(k)\,,
\quad a(k) \geq 0\,,
\end{equation}
whose partial sums are given by (\ref{defsn}), into an alternating
series $\sum_{j=0}^{\infty} (-1)^j {\bf A}_j$. After the first step
of the transformation, the limit of the input series is recovered
according to
\begin{equation}
\label{VWSerTran}
\sum_{k=0}^{\infty} \, a (k) =
\sum_{j=0}^{\infty} \, (-1)^j \, {\bf A}_j \,.
\end{equation}
The quantities ${\bf A}_j$ are defined according to 
\begin{equation}
\label{A2B}
{\bf A}_j = \sum_{k=0}^{\infty} \, {\bf b}_{k}^{(j)} \, ,
\end{equation}
where
\begin{equation}
\label{B2a}
{\bf b}_{k}^{(j)} = 2^k \, a(2^k\,(j+1)-1) \, .
\end{equation}
The ${\bf A}_j$ are referred to as the {\em condensed
series}~\cite{JeMoSoWe1999}, and the series $\sum_{j=0}^{\infty}
(-1)^j {\bf A}_j$ is referred to as the {\em transformed alternating
series}, or alternatively as the {\em Van Wijngaarden transformed
series}.

The construction of the condensed series reminds one of Cauchy's
condensation theorem (see e.g.~p.~28 of Ref.~\cite{Bw1991} or p.~121
of Ref.~\cite{Kn1964}). Given a nonalternating series
$\sum_{k=0}^{\infty} a (k)$ with terms that satisfy $\vert a(k+1)
\vert < \vert a(k) \vert$, Cauchy's condensation theorem states that
$\sum_{k=0}^{\infty} a (k)$ converges if and only if the first
condensed series ${\bf A}_0$ defined according to Eq.~(\ref{A2B})
converges.

The summation over $k$ in Eq.~(\ref{A2B}) does not pose numerical
problems. Specifically, it can be easily shown in many cases of
practical importance that the convergence of $\sum_{k=0}^{\infty} \,
{\bf b}_{k}^{(j)}$ (in $k$) is linear even if the convergence of
$\sum_{k=0}^{\infty} a(k)$ is only logarithmic. We will illustrate
this statement by way of two examples. Example 1: a logarithmically
convergent input series whose terms behave asymptotically as $a (k)
\sim k^{-1 - \epsilon}$ with $\epsilon > 0$. In this case, the
partial sums
\begin{equation}
{\bf A}^{(n)}_j = \sum_{k=0}^{n} \, {\bf b}_{k}^{(j)}
\end{equation}
converge linearly with
\begin{equation}
\lim_{n \to \infty} \frac{{\bf A}^{(n+1)}_j - {\bf A}_j} 
  {{\bf A}^{(n)}_j - {\bf A}_j} \; = \; 
\frac{1}{2^\epsilon \, (j+1)^{1+\epsilon}} < 1 \,, \qquad
a(k) \sim k^{- 1 - \epsilon}\,, \quad k \to \infty.
\end{equation}
Example 2: a series with $a (k) \sim k^{\beta} r^k$ where $0 < r <
1$ and $\beta$ real. Here, we have $\rho = r < 1$, and the series is
(formally) linearly convergent. However, slow convergence may result
if $\rho$ is close to one. In this case, the condensed series are
very rapidly convergent,  
\begin{equation}
\label{criterion}
\lim_{n \to \infty} \frac {{\bf A}^{(n+1)}_j - {\bf A}_j} 
  {{\bf A}^{(n)}_j - {\bf A}_j} \; = \; 0 \,, \qquad
a(k) \sim k^{\beta} r^k \,, \quad k \to \infty.
\end{equation}
Therefore, when summing over $k$ in evaluating the condensed series
according to Eq.~(\ref{A2B}), it is in many cases sufficient to 
to evaluate the condensed series by adding the terms successively,
and no further acceleration of the convergence is required.

As shown in~\cite{Da1969,JeMoSoWe1999}, the condensation
transformation defined according to
Eqs.~(\ref{VWSerTran})--(\ref{B2a}) is essentially a reordering of
the terms of the input series $\sum_{k=0}^{\infty} a(k)$.
Furthermore, Daniel was able to show (see the Corollary on p.~92 of
Ref.~\cite{Da1969}), that for nonalternating convergent series whose
terms decrease in magnitude ($\vert a(k) \vert > \vert a(k+1)
\vert$), the equality (\ref{VWSerTran}) holds. This formally
justifies the correctness of the condensation transformation defined
according to Eqs.~(\ref{VWSerTran}) -- (\ref{B2a}). 

Note that the property, originally derived in~\cite{Da1969},
\begin{equation}
\label{half}
{\bf A}_{2\,j-1} = \frac{1}{2} \,
\left( {\bf A}_{j-1} - a_{j-1} \right)\,, \qquad
\mbox{($j = 1,2,\dots$)}\,,
\end{equation}
facilitates the numerical evaluation of a
set of condensed series, by reducing the evaluation of condensed
series of odd index to a trivial computation.
Within the program, we use this relation in the form
\begin{equation}
\label{compute}
{\bf A}_{i+1} = \frac{1}{2}\,({\bf A}_{i/2} - a_{i/2}) \qquad
\mbox{(even $i$)}
\end{equation}
in order to compute condensed series $A_m$ of odd index $m=i+1$
where $i$ is even.

In the second step of the CNCT, the convergence of the Van
Wijngaarden transformed series $\sum_{j=0}^{\infty} \, (-1)^j \,
{\bf A}_j$ on the right-hand side of Eq.~(\ref{VWSerTran}) is
accelerated by a suitable nonlinear sequence transformation. We
start from the partial sums 
\begin{equation}
\label{PSumS}
{\bf S}_n \; = \; \sum_{j=0}^{n} \, (-1)^j \, {\bf A}_j
\end{equation}
of the Van Wijngaarden transformed series. Formal operations given
by operators acting on the space of sequences are of considerable
importance for the construction and theory of sequence
transformations. We define the (forward) difference operator $\Delta$
\begin{equation}
\label{defOpDelta}
\Delta \; : \; 
\{s_n\}_{k=0}^{\infty} \to \{(\Delta s)_n\}_{n=0}^{\infty} \; , \;
(\Delta s)_n \equiv s_{n+1} - s_n\,.
\end{equation}
The $n$th element of the sequence $\{(\Delta s)_n\}_{n=0}^{\infty}$
is given by $(\Delta s)_n = s_{n+1} - s_n$. This way of writing the
$n$th element of the {\em transformed} sequence stresses the fact
that $(\Delta s)_n$ is to be regarded as an element of some new
sequence. However, the brackets are often left out, and we write
\begin{equation}
\Delta s_n \equiv s_{n+1} - s_n = (\Delta s)_n\,.
\end{equation}
A well-known important relation is
\begin{equation}
\label{DeltaToTheK}
(\Delta^k s)_n = (-1)^k \,
\sum_{j = 0}^{k} (-1)^j \, 
\left( \begin{array}{c} k \\ j \end{array} \right) \, 
s_{n+j}\,.
\end{equation}
As before we identify
\begin{equation}
\Delta^k s_n \equiv (\Delta^k s)_n\,.
\end{equation}
The $k$th power of the difference operator plays a crucial rule in
constructing sequence transformations~[see e.g.~\cite{We1989}].
A very important class of sequence transformations 
(``Levin-type transformations''), which are characterized
by the use of explicit remainder estimates $\omega_n$, is constructed
according to the following 
``prescription''~\cite{Le1973,Si1981,We1989,Ho2000}
\begin{eqnarray}
\label{SidTr}
\lefteqn{
{\cal S}_{k}^{(n)} (\beta, s_n, \omega_n) \; = \; \frac
{ \Delta^k \, \{ (n + \beta)_{k-1} \; s_n / \omega_n\} }
{ \Delta^k \, \{ (n + \beta)_{k-1}  / \omega_n \} }} \nonumber \\[2ex]
& \; = \; \frac
{\displaystyle
\sum_{j=0}^{k} \; ( - 1)^{j} \; 
\left( \begin{array}{c} k \\ j \end{array} \right) \; 
\frac {(\beta + n +j )_{k-1}} {(\beta + n + k )_{k-1}} \;
\frac {s_{n+j}} {\omega_{n+j}} }
{\displaystyle
\sum_{j=0}^{k} \; ( - 1)^{j} \; 
\left( \begin{array}{c} k \\ j \end{array} \right) \; 
\frac {(\beta + n +j )_{k-1}} {(\beta + n + k )_{k-1}} \;
\frac {1} {\omega_{n+j}} } \, .
\end{eqnarray}
The quantities appearing in the definition of ${\cal S}_{k}^{(n)}$
have the following interpretation:
\begin{itemize}
\item $s_n$: elements of a sequence whose convergence is to
be accelerated.
\item $\omega_n$: remainder estimates, i.e.~estimates for the
truncation error $\omega_n \approx r_n$ defined in (\ref{defrn}).
Various sequence transformations based on 
Eq.~(\ref{SidTr}), which differ only in the choice of
$\omega_n$, have been discussed in detail in~\cite{We1989}. For all
these sequence transformations, $\omega_n$ can be
calculated on the basis of the $s_n$. E.g., for
the delta transformation defined below in Eq.~(\ref{defdelta}), 
which is exclusively studied in the following, the calculation
of $\omega_n = \Delta s_n = s_{n+1} - s_n$ requires as input 
the sequence elements $s_n$ and $s_{n+1}$. 
\item $\beta$: a numerical shift parameter, which may be adjusted or
optimized for a particular application, e.g.~by establishing
asymptotic relations between the $\omega_n$ and the $r_n$, or by
numerical experimentation. In practice, it is observed that any
choice other than $\beta = 1$ does not lead to an appreciable
improvement of the rate of convergence of the transforms ${\cal
S}_{k}^{(n)}$, and therefore $\beta$ is almost exclusively set to
unity in the literature (see also the numerical examples
in~\cite{We1989,JeMoSoWe1999}). In order to allow for the case $n=0$ in
Eq.~(\ref{SidTr}), we must demand that $\beta \neq 1$, and 
in order to ensure a regular behaviour of the remainder estimates
as a function of the index $n+j$, we should choose $\beta > 0$.
\item $k$: order the transformation. The forward difference operator
$\Delta$ acts on the index $n$. The $k$th power of that difference
operator enters into the numerator and into the denominator of
Eq.~(\ref{SidTr}).
\item $n$: ``initial element'' of the transform. The computation of
the transform ${\cal S}_{k}^{(n)}$ requires as input the $k+1$
elements $\{ s_n,\dots, s_{n+k} \}$ of the sequence $\{
s_n \}_{n=0}^{\infty}$, together with the $k+1$
remainder estimates $\{ \omega_n,\dots, \omega_{n+k} \}$. The element
of lowest index required for the calculation of ${\cal S}_{k}^{(n)}$
is $s_n$, and therefore $n$ can be interpreted as the ``initial
element'' in the calculation of the transform ${\cal S}_{k}^{(n)}$.
\end{itemize}
At present, no mathematical proof exists with regard to the
convergence of the transforms ${\cal S}_{k}^{(n)}$,
although the construction of sequence transformations of the type
(\ref{SidTr}) can be motivated
heuristically  (by expanding the  
truncation error in a factorial series) 
and is supported by overwhelming numerical
evidence~\cite{We1989,We1996c,Je1999,JeMoSoWe1999}. 

The notation ${\cal S}_{k}^{(n)} (\beta, s_n, \omega_n)$ is in need
of a certain further explanation. As explained above, $n$ represents
the initial element of the sequence $\{ s_n \}_{n=0}^{\infty}$ used
in the evaluation of the transform. The specification of $s_n$ 
as the second argument of ${\cal S}_{k}^{(n)} (\beta, s_n, \omega_n)$
is therefore redundant as far as the index $n$ is concerned, because
$n$ already appears as an upper index in ${\cal S}_{k}^{(n)}$. The
arguments $s_n$ (and $\omega_n$) are to be interpreted as follows:
they rather specify the initial element $s_n$ and the initial
remainder estimate $\omega_n$ which are required for the
calculation of the right-hand side of (\ref{SidTr}). 
Of course, the evaluation of ${\cal S}_{k}^{(n)} (\beta, s_n, \omega_n)$ 
requires in total the $k+1$ sequence
elements $\{ s_n,\dots, s_{n+k} \}$ and the $k+1$
remainder estimates $\{ \omega_n, \dots, \omega_{n+k} \}$
as input. The knowledge of the ``post-initial'' sequence
elements $\{ s_{n+1},\dots, s_{n+k} \}$ and remainder estimates
$\{ \omega_{n+1},\dots, \omega_{n+k} \}$
is by assumption guaranteed in writing the definition
(\ref{SidTr}). A somewhat more clumsy notation for
${\cal S}_{k}^{(n)} (\beta,  s_n, \omega_n)$ 
which avoids all redundancy would read
\[
{\cal S}_{k}^{(n)} (\beta, \{s\}_{j=0}^\infty, \{\omega\}_{j=0}^\infty )
\qquad \mbox{(alternative notation not used in this article)}.
\]

The delta transformation~\cite{Si1981,We1989} is a special case of 
the transformation~(\ref{SidTr}). It is constructed according to
\begin{equation}
\label{defdelta}
\delta_{k}^{(n)} \left( \beta, s_n \right) \; = \;
{\cal S}_{k}^{(n)} \left(\beta, s_n, \Delta s_n \right) \, .
\end{equation}
We use as input data for the delta transformation the partial sums
(\ref{PSumS}) of the Van Wijngaarden transformed series. This
corresponds to the ``replacements''
\begin{equation}
\label{repl}
s_n \to {\bf S}_n \,, \qquad 
\Delta s_n \to \Delta {\bf S}_n = (-1)^{n+1} \, {\bf A}_{n+1}\,,
\end{equation}
and here, we always use ${\bf S}_0$ as the initial element for the 
second transformation (see however the discussion on p.~203
of~\cite{Je1999} where it is indicated that some terms of a 
given input series may have to be
skipped in practical applications where the first terms of 
the input series
display a highly irregular behaviour). This leads to the CNC transforms
\begin{equation}
\label{defcnc}
{\cal T}_{\rm CNC}(n) = \delta_{n}^{(0)} \left( 1, {\bf S}_0 \right)\,,
\end{equation}
which require as input the elements
$\{ {\bf S}_0, \dots, {\bf S}_{n}, {\bf S}_{n+1} \}$
of the Van Wijngaarden transformed series. 

A final word on sequence transformations: It is observed that
sequence transformations constructed according to (\ref{SidTr})
often lead to much better numerical results than better known
convergence accelerators such as the Euler transformation, Wynn's
epsilon and rho algorithms, or Aitken's $\Delta^2$ process
(all of these algorithms are described in~\cite{We1989}). In view
of the partially unsatisfactory situation concerning the
availability of mathematical proof, we have carried out
high-precision studies of the convergence of sequence
transformations of the type (\ref{SidTr}), in part using
multi-precision arithmetic. In the absence of a rigorous proof, one
of the concerns which may be raised against the convergence of the
transforms (\ref{SidTr}) is a plausible asymptotic nature of the
sequence of transforms, resulting in ultimate divergence in
higher order. In our numerical experiments, we found no
indication of unfavorable
asymptotic behaviour. The rate of convergence of the
transforms remained constant, and we observed apparent convergence
up to 430 decimal figures (see Sec.~\ref{sec_math} below). In the
absence of rigorous proof and the presence of considerable numerical
evidence for the computational usefulness of sequence
transformations of the type (\ref{SidTr}), it is conceivable that
``mathematical experiments'' may give some hints at the theoretical
soundness of the concepts involved.

%
%
\section{Applications to Special Functions used in Statistics}
\label{sec_stat}

Several slowly convergent series of the type considered in this
article define special functions that have important applications in
statistics. These functions are Riemann's zeta, 
the generalized zeta, Jonqui\`{e}re's
function, and Lerch's transcendent (see~\cite{Ba1953vol1} for
definitions and some properties of these functions).  Discrete
distributions are related to these functions by associating
probability mass functions (p.m.f.s) with the terms of 
the infinite series that define the functions.
In turn, the normalization constants of the p.m.f.s are associated
with the sums of the series.

Consider a discrete distribution with the p.m.f. 
\begin{equation} 
\label{discrete}
{\rm Pr}[X=n] = c\,p(n;\theta)\,,
\end{equation}
and with support of all nonnegative integers $n=0,1,\dots$, where
$\theta$ is a (vector of) parameter(s) and $c$ is a constant.  
The probability over the set of all
outcomes is unity:
\begin{equation} 
\label{proper}
\sum_{n=0}^{\infty} {\rm Pr}[X=n] = 1\,,
\end{equation}
from which we obtain the normalization constant
\begin{equation} 
\label{normc}
c^{-1} = \sum_{n=0}^{\infty} p(n;\theta)
\end{equation}
(in statistics, the random variable is usually denoted by the symbol $x$,
whether continuous or discrete, but we will use the symbol $n$
here, which is the more common notation for the running index
in the theory of special functions).
Due to its very property as a normalization constant,
we expect the sum in Eq.~(\ref{normc}) to exist and to converge
to a finite value. In the following,
we consider several related discrete distributions.
The Zipf distribution has the p.m.f.
\begin{equation}
\label{zipf}
p_n = c\,\frac1{n^s}
\end{equation}
with support of all positive integers $n=1,2,\dots$.  The
normalization constant~(\ref{normc}) is then
\begin{equation} 
\label{zeta}
c^{-1} = \zeta(s) = \sum_{n=1}^\infty \frac{1}{n^s}\,,
\end{equation}
where $s>1$. This is immediately recognized as the Dirichlet series
for Riemann's zeta function [Eq.~(1) on p.~32 of
Ref.~\cite{Ba1953vol1}].  Basic properties of the Zipf distribution
can be found in~\cite{JoKoKe1992}. 

The Zipf--Mandelbrot distribution~\cite{Ma1983} is the generalization 
of the Zipf
distribution that adds a constant $v$ to the ranks $n$ and has the p.m.f. 
\begin{equation}
\label{zipfman}
p_n=c\,\frac1{(n+v)^s}\,,
\end{equation}
with the support of all nonnegative integers $n=0,1,\dots$. The
normalization constant~(\ref{normc}) is then
\begin{equation} 
\label{genzeta}
c^{-1} = \zeta(s,v) = \sum_{n=0}^\infty \frac1{(n+v)^s}\,,
\end{equation}
where $v \neq 0,-1,\dots$. The sum~(\ref{genzeta}) is recognized as
the generalized zeta function [Eq.~(1) on p.~24 of
Ref.~\cite{Ba1953vol1}].

Further, the p.m.f.
\begin{equation}
\label{good}
p_n=c\,\frac{z^n}{n^s}
\end{equation}
with the support of all positive integers $n=1,2,\dots$ defines the
Good distribution~\cite{Go1953}. The normalization
constant~(\ref{normc}) is then
\begin{equation}
\label{jonquiere}
c^{-1} = F(z,s) = \sum_{n=1}^\infty \frac{z^n}{n^s}\,,
\end{equation}
where $\vert z \vert <1$, which is recognized as Jonqui\`{e}re's
function [Eq.~(14) on p.~30 of Ref.~\cite{Ba1953vol1}; see 
also~\cite{Tr1945}].  Kulasekera and Tonkyn~\cite{KuTo1992} 
derived moments of the Good distribution and proposed parameter 
estimators (for parameter estimation see also~\cite{DoLu1997}).

Now we discuss briefly several applications of the above
distributions.  While Zipf and Zipf--Mandelbrot distributions have
been traditionally used in linguistics and information science, as
well as in other disciplines~\cite{Zi1949,Ma1983}, they also appear
as canonical distributions in the generalized statistical-mechanics
ensemble theory.  A variety of physical systems e.g. with long-range
interactions or fractal boundary conditions are intractable within
the classical Boltzmann--Gibbs statistical mechanics, while they can
be successfully handled within the generalization proposed by
Tsallis~\cite{Ts1988}.  In particular, the variation of the
generalized entropy with appropriate constraints yields power-law
equilibrium distributions of the Zipf--Mandelbrot
type~\cite{TsMePl1998}. The convolution of these distributions was
shown to converge to the limiting L\'{e}vy
distribution~\cite{AbRa2000}, according to the L\'{e}vy-Gnedenko
generalized central-limit theorem~\cite{GnKo1968}.  Furthermore,
the interrelationship between microscopic statistical and macroscopic
thermal equilibria for such power-law systems was shown to be a
natural extension of the corresponding Boltzmann--Gibbs
theory~\cite{AbRa2001}.

The Good distribution has been used in linguistics to model word
frequencies~\cite{Si1955,He1958}. It was originally formulated by
Good~\cite{Go1953} by introducing a ``convergence factor'' $z^n$ in
Zipf's p.m.f.~(\ref{zipf}) in order to make the first moment of
the distribution finite for $1 < s \leq 2$.  The p.m.f.~(\ref{good})
as an approximation for the upper tail of the Yule
distribution~\cite{Yu1924} was mentioned by Simon~\cite{Si1955}. 
The particular case of the Good distribution with $s=1$, known as 
Fisher's logarithmic law~\cite{FiCoWi1943}, has been obtained as an
equilibrium solution of certain stochastic processes by Champernowne
in his model of the distribution of incomes~\cite{Ch1953}, and by
Darwin in his model of the frequencies of species in evolving
populations~\cite{Da1953}. A logarithmic law for the frequencies of
species has also been obtained as a limiting form of the negative
binomial distribution~\cite{Wi1944}. The Good distribution was also
suggested for use in ecology for the reasons of the flexibility of
its hazard function~\cite{KuTo1992}.

In biophysics, one is interested in establishing the statistical
structure of various sequences, e.g. those of DNA, RNA, or proteins. 
Indeed, an interest in statistical compositions of DNA and protein sequences
originated shortly after the discovery that the nucleotide sequences
of RNA uniquely determines the amino-acid sequence of the proteins. 
An early analysis was performed by Gamow and
Y\v{c}as~\cite{GaYc1955}, who observed that the distributions of
relative abundances of amino-acids and nucleotides were nonrandom,
i.e. deviated from models assuming a uniform distribution of different
words in any given text.  Further studies established the utility of
relative abundances of words (short oligonucleotides) in a comparative
analysis of genome sequences~\cite{KaBu1995,BlCaKa1996}.  Furthermore, a
linguistically motivated
``Zipf'' analysis has been applied to a number of DNA
sequences in an attempt to demonstrate that noncoding sequences
resemble a ``natural'' language more than the coding
ones~\cite{Maetal1994}. However, it was noted that the differences
between the fitted Zipf exponents $s$ in the p.m.f.~(\ref{zipf}) are
not indicative of the closeness to a ``natural'' language per
se~\cite{KoMa1995}, as ``Zipf--like'' (i.e., power-law) frequency
distributions arise in a wide variety of situations.  Moreover, the
obtained fits do not apply in the whole range of word ranks (i.e.,
slopes of empirical plots in log-log coordinates are not constant).  In
contrast, the Good distribution was found to fit the data in a
wider range of ranks~\cite{MaKo1996}, even though it was not
possible to determine the statistically significant differences
between the fits to coding and noncoding sequences, due to large
sampling errors of the fitted parameters.  These studies reinforce an
early comment by Simon~\cite{Si1955} that the possible cause for the
wide applicability of the ``power-law'' Zipf-related distributions
is perhaps the similarity in the probability mechanisms underlying
such diverse phenomena as distributions of words in DNA sequences and
in ``natural'' languages, frequencies of citations of scientific
articles and of species abundances in populations.  Interestingly, the
possibility of such a unification might occur within the framework of
the generalized thermostatistics, which is exemplified by obtaining
the Zipf--Mandelbrot law in binary sequences within the generalized
thermodynamics of Tsallis~\cite{De1997}.

A formal generalization of the Zipf, the Zipf--Mandelbrot, and the
Good distributions is possible upon realizing that Riemann's zeta,
the generalized zeta, and Jonqui\`{e}re's functions constitute
special cases of Lerch's transcendent which is defined by the
following series
\begin{equation}
\label{lerch}
\Phi(z,s,v) = \sum_{n=0}^\infty \frac{z^n}{(n+v)^s}\,,
\end{equation}
where $|z|<1$ and $v \neq 0,-1,\dots$ [Eq.~(1.11.1) on p.~27 of
Ref.~\cite{Ba1953vol1}]. The relations between $\Phi(z,s,v)$ and
$\zeta(s)$, $\zeta(s,v)$, and $F(z,s)$ can be obtained by letting $z$
and/or $v$ to 0 and/or 1, and making use of the following functional
relation [Eq.~(1.11.2) on p.~27 of Ref.~\cite{Ba1953vol1}]
\begin{equation}
\label{lerchrel}
\Phi(z,s,v) = z^m \, \Phi(z,s,m+v) + \sum_{n=0}^{m-1} \frac{z^n}{(n+v)^s}\,.
\end{equation}
In particular, Eq.~(\ref{lerchrel}) with $m=1$ can be used
to rewrite sums that start from $n=1$ rather than from $n=0$
in terms of the $\Phi$ function:
\begin{equation}
\sum_{n=1}^\infty \frac{z^n}{(n+v)^s} = z \, \Phi(z,s,1+v)
\end{equation}

\begin{table}
\begin{center}
\begin{minipage}{14cm}
\begin{center}
\begin{tabular}{llll}
\hline
\hline\\[0.3ex]
Distribution & Support & p.m.f. & normalization \\
(related function) & & & constant $c^{-1}$ \\[2ex]
\hline\\[0.3ex]
Zipf & 
$1,2,\dots$ & 
$c\,n^{-s}$ & 
$\zeta(s) = \Phi(1,s,1)$\\
(Riemann's zeta) & & & \\[2ex]
Zipf--Mandelbrot & 
$0,1,\dots$ & 
$c\,(n+v)^{-s}$ & 
$\zeta(s,v) = \Phi(1,s,v)$\\
(generalized zeta) & & & \\[2ex]
Good & 
$1,2,\dots$ & 
$c\,z^n\,n^{-s}$ & 
$F(z,s) = z \, \Phi(z,s,1)$\\
(Jonqui\`{e}re's function) & & & \\[2ex]
Lerch & 
$0,1,\dots$ & 
$c\,z^n\,(n+v)^{-s}$ & 
$\Phi(z,s,v)$\\
(Lerch's transcendent) & & & \\[0.3ex]
\hline
\hline
\end{tabular}
\end{center}
\caption{Relationships between distributions defined by 
Riemann's zeta, generalized zeta, Jonqui\`{e}re's and Lerch's functions.}
\label{speccases}
\end{minipage}
\end{center}
\end{table}

In Tab.~\ref{speccases}, we have summarized the relations between
the different statistical distributions and their defining special
functions. The three distributions, Zipf (zeta), Zipf--Mandelbrot
(generalized zeta), and Good (Jonqui\`{e}re), are special cases of
the Lerch distribution, and their properties can be expressed in
terms of Lerch's transcendent with special values of the parameters. 
(See~\cite{ZoAl1995} for the discussion of relations between these 
distributions and the zero-truncated Lerch distribution.)

Here we consider several functions of the related Lerch distributions,
from which other properties (e.g., moments) can be easily derived
(see~\cite{AkSa2001} for details on the Lerch distribution).

The cumulative distribution function (c.d.f.) $F(n;\theta)$ of 
the Lerch distribution with the support
$n=0,1,\dots$ is given by the following equation
\begin{equation} 
\label{lerchdf}
F(n;\theta) = {\rm Pr}[X \leq n] = 1 - z^{n+1} \,
\frac{\Phi(z,s,v+n+1)}{\Phi(z,s,v)}\,.
\end{equation}
The c.d.f.s
for the Zipf, the Zipf--Mandelbrot, and the Good distributions are given by
\begin{equation} 
\label{specdf}
1-\frac{\Phi(1,s,n+1)}{\Phi(1,s,1)},\qquad  
1-\frac{\Phi(1,s,v+n+1)}{\Phi(1,s,v)},\qquad 
1-\frac{z^n \Phi(z,s,n+1)}{\Phi(z,s,1)}\,,
\end{equation}
respectively.  The survival and hazard functions are obtained from
the c.d.f.s~(\ref{lerchdf}) and~(\ref{specdf}) as
$S(n;\theta)=1-F(n;\theta)$ and
$h(n;\theta)=p(n;\theta)/(1-F(n;\theta))$, respectively.

The probability generating function (p.g.f.) $G(y;\theta)$ for the 
Lerch distribution is given by
\begin{equation} 
\label{lerchpgf}
G(y;\theta) = \sum_{n=0}^{\infty} y^n \, {\rm Pr}[X=n] = 
\frac{\Phi(y\,z,s,v)}{\Phi(z,s,v)}\,,
\end{equation}
where we assume $\vert y \vert \leq 1$. The p.g.f.s for the Zipf,
Zipf--Mandelbrot and Good distributions are then given by
\begin{equation} 
\label{specpgf}
\frac{y \, \Phi(y,s,1)}{\Phi(1,s,1)},
\qquad 
\frac{\Phi(y,s,v)}{\Phi(1,s,v)},
\qquad 
\frac{y \, \Phi(y\,z,s,1)}{\Phi(z,s,1)}\,,
\end{equation}
respectively.
The characteristic functions (c.f.'s) and moment generating functions
(m.g.f.s) can be obtained from the p.g.f.s~(\ref{lerchpgf})
and~(\ref{specpgf}) by letting $y=e^{it}$ and $y=e^t$, respectively.
The moments of the Lerch distributions are then obtained as
coefficients of the Taylor series of the m.g.f.s about $t=0$.

We note that so far the support of the distributions consisted of
all nonnegative integers, except possibly 0, whereas in applications
the distributions are often used in singly or doubly truncated
forms. Truncation of a distribution changes the way the normalizing
constant $c^{-1}$ is calculated. For the truncated case of the Lerch
distribution with support $n\in[a,b]$, where $a \geq 0$ and $a \leq
b \leq \infty$, the c.d.f. is now~\cite{AkSa2001}
\begin{equation} 
\label{lerchdftr}
F(n;\theta) = c \sum_{k=a}^n \frac{z^k}{(k+v)^s} = 
\frac{z^a \, \Phi(z,s,v+a) - z^{n+1} \, \Phi(z,s,v+n+1)}
  {z^a \, \Phi(z,s,v+a) - z^{b+1} \, \Phi(z,s,v+b+1)}\,.
\end{equation}
Likewise, the p.g.f. for the Lerch distribution is
\begin{equation}
\label{lerchpgftr}
G(y;\theta) = 
\frac{(y\,z)^a \, \Phi(y\,z,s,v+a) - (y\,z)^{b+1} \, \Phi(yz,s,v+b+1)}
  {z^a \, \Phi(z,s,v+a) - z^{b+1} \, \Phi(z,s,v+b+1)}\,.
\end{equation}
For truncated Zipf, Zipf--Mandelbrot, and Good distributions, the c.d.f.s
can be obtained from Eq.~(\ref{lerchdftr}), and the p.g.f.s can be
obtained from Eq.~(\ref{lerchpgftr}), by substituting $z=1$ and
$v=0$ (for Zipf), $z=1$ (for Zipf--Mandelbrot) and $v=0$ (for Good).
Note that these ``short-cut'' substitutions are possible because
truncated distributions have the same support (by definition). 
We cannot use these
substitutions to obtain Eq.~(\ref{specdf}) and~(\ref{specpgf})
because in that case the distributions have different support,
depending on whether zero is excluded or included (see Tab.~\ref{speccases}).

Calculations with the Lerch distributions using
Eqs.~(\ref{lerchdf}) -- (\ref{lerchpgftr}) involve an evaluation of
Lerch's transcendent.  Other commonly needed manipulations with the
Lerch distributions, e.g. random number generation and parameter
estimation, can be likewise reduced to calculating the Lerch
transcendent~\cite{AkSa2001}.

We have implemented the calculation of Lerch's transcendent $\Phi(z,s,v)$ 
defined in Eq.~(\ref{lerch}) using the 
CNC transformation [Eqs.~(\ref{VWSerTran}) -- (\ref{B2a})] in 
Mathematica. The program is presented and discussed in Appendix~\ref{app}.
For a {\em provisional} implementation in C
as well as a Mathematica
package for calculations with the Lerch distribution, 
we refer to~\cite{AkJeURL}.

%
%
\section{An Application in Experimental Mathematics}
\label{sec_math}

We begin this section by 
quoting~\cite{Ba1994tech}: ``In
April 1993, Enrico Au--Yeung, an undergraduate at the University of
Waterloo, brought to the attention of the author's [David Bailey's]
colleague Jonathan Borwein the curious fact that 
\begin{equation}
\label{SumBailey} 
\sum_{k=1}^\infty \;  \left( 1 + \frac{1}{2} + \cdots +
\frac{1}{k} \right)^2 \, k^{-2} \;=\; 4.59987\cdots \; \approx \;  
\frac{17}{4} \, \zeta(4) \;=\;  \frac{17 \pi^2}{360} \,, 
\end{equation} 
based on a computation of 500,000 terms. Borwein's reaction was to
compute the value of this constant to a higher level of precision in order
to dispel this conjecture.  Surprisingly, his computation to 30 digits
affirmed it. [David Bailey] then computed this constant to 100 decimal
digits, and the above equality was still affirmed.''
Many formulas similar to (\ref{SumBailey}) have subsequently been established
by rigorous proof~\cite{BaBoGi1994}. 

With the help of a multiprecision system~\cite{Ba1990tech,Ba1993,Ba1994tech} 
and the CNCT, we have verified (\ref{SumBailey}) ``experimentally''
to a couple of hundred decimals~\cite{Be1999master}. The calculation will be
sketched in the following. Using the definition
\begin{equation}
{\bar b}(k) \;=\; 
\left( 1 + \frac{1}{2} + \ldots + \frac{1}{k} \right)^2 \, k^{-2}\,,
\end{equation}
we rewrite (\ref{SumBailey}) as follows,
\begin{eqnarray}
\sum_{k=0}^\infty \; {\bar b}(k)
& = & 
\sum_{k=0}^\infty \; 
\left( \sum_{j=0}^{k} \; \frac{1}{j+1} \right)^2 \, (k+1)^{-2}
\nonumber \\[2ex]
& = &
\sum_{k=0}^\infty \; \left( \frac{\psi(k+2) + \gamma}{k+1} \right)^2 \,,
\end{eqnarray}
where $\gamma$ is the {\em Euler-Mascheroni} constant
$\gamma = 0.577~215\ldots$, and $\psi(z)$ is the 
logarithmic derivative of the 
Gamma function~\cite{Ol1974,Ba1953vol1}\,,
\begin{equation}
\psi(z) \;=\; \frac{\mathrm d}{{\mathrm d}z} \ln \Gamma(z)\,.
\end{equation}
With the help of the relation
\begin{equation}
\sum_{k=0}^\infty \; \frac{\psi(k+2)}{(k+1)^2} =
  2 \, \zeta(3) - \gamma \, \zeta(2)\,,
\end{equation}
Eq.~(\ref{SumBailey}) can be rewritten as
\begin{equation}
\label{SimplifySumBailey}
\sum_{k=0}^\infty \; \left( \frac{\psi(k+2)}{k+1} \right)^2 =
\frac{17}{4} \, \zeta(4) - 4 \gamma\,\zeta(3) + \gamma^2 \, \zeta(2)\,.
\end{equation}
We proceed to calculate numerically, to high
precision, the infinite sum
\begin{equation}
\label{SumTodo}
\sum_{k=0}^\infty \; {\bar a}(k) \,, \qquad
{\bar a}(k) = \left( \frac{\psi(k+2)}{k+1} \right)^2\,,
\end{equation}
using the CNC transformation.

In order to establish the rate of convergence
of (\ref{SumTodo}), we investigate the asymptotic
behaviour of the ${\bar a}(k)$ as $k \to \infty$.
The logarithm of the Gamma function can be expanded into an 
asymptotic series [see Eq.~(4.03) on p.~294 of~\cite{Ol1974}]:
\begin{equation}
\label{LogGamma2} 
\ln \Gamma(z) \;=\; \left(z - \frac{1}{2} \right) \, \ln z
\;-\; z \;+\; 
\frac{1}{2} \, \ln (2 \pi) \;+\; 
\sum_{s=1}^{m-1} \; 
\frac{{\cal B}_{2s}}{2s(2s-1) z^{2s-1} } \;+\; R_m(z) \,, 
\end{equation}
where
\begin{equation}
R_m(z) \;=\; 
\int_0^\infty \frac{ {\mathcal B}_{2m} -
  {\mathcal B}_{2m}(x - \Ent{x})} 
    {2 m (x + z)^{2m}} {\mathrm d}x \,. 
\end{equation}
Here, $\Ent{x}$ is the integral part of
$x$, i.e., the largest integer $m$ satisfying $m \le x$, 
$B_k (x)$ is a Bernoulli polynomial defined by the 
generating function 
[see Eq.~(1.06) on p.~281 of Ref.~\cite{Ol1974}]:
\begin{equation}
\frac {t \exp(x t)}{\exp(t) - 1} \; = \;
\sum_{m=0}^{\infty} \, B_m(x) \, \frac {t^m} {m!} \, ,
\qquad \vert t \vert < 2 \pi \, ,
\end{equation}
and
\begin{equation}
B_m \; = \; B_m (0)
\label{BernNum}
\end{equation}
is a Bernoulli number (p.\ 281 of Ref.~\cite{Ol1974}).
The following asymptotic relation for $\psi(z)$ follows:
\begin{equation}
\label{PolyAsym}
\psi(z) = \ln z \;-\; 
\frac{1}{2z} \;-\; \sum_{s=1}^{m-1} \; 
\frac{{\mathcal B}_{2s}}{2s z^{2s} } \;+\; 
{\mathcal O}\left(\frac{1}{z^{2m}} \right)\,.
\end{equation}
The leading asymptotics of
the remainder of the 
sum (\ref{SumTodo}) after adding $N-1$ terms can thus be derived easily.
We have for large $k$,
\begin{equation} 
\label{asymp_a}
{\bar a}(k) \;  \sim \; 
\frac{\ln(k + 2)^2}{(k+1)^2} -
\frac{\ln(k + 2)}{(k + 1)^2 \, (k + 2)} -
\frac{\ln(k + 2)}{6\,(k + 1)^2 \, (k + 2)^2}
+ {\mathcal O}\left(\frac{1}{k^4} \right)  \,,
\qquad k \to \infty\,.
\end{equation}
Based on these formulas, 
the remainder of the sum (\ref{SumTodo}), for large $N$, 
can be written as
\begin{equation}
\label{emestimate}
\sum_{k=N}^\infty \; {\bar a}(k) \; \sim \;
\frac{\ln^2 N}{N} \;+\; \frac{\ln N}{N} \;+\; \frac{1}{N} \;+\;
{\mathcal O}\left( \frac{ \ln^2 N}{N^2} \right) \,. 
\end{equation}
Here, the Euler-Maclaurin formula [Eqs.~(2.01)
and~(2.02) on p.~285 of Ref.~\cite{Ol1974}] has been used in order
to convert the sum over the ${\bar a}(k)$ in the
asymptotic regime of large $k$ [see Eq.~(\ref{asymp_a})]
into an integral plus correction terms.
In order to calculate (\ref{SimplifySumBailey}) to an accuracy of
200 decimals, Eq.~(\ref{emestimate}) says
that we would be required to add
on the order of $10^{205}$ terms. Without the use 
of convergence acceleration methods, this would represent
a formidable computational task. 

Using the CNCT, it is easy to
calculate the sum (\ref{SimplifySumBailey}) to 200 digits, 
based on multiprecision
arithmetic~\cite{Ba1990tech,Ba1993,Ba1994tech} 
and a Linux personal computer, within a few
hours. We obtain for the 246th and the 247th CNC transform defined
according to Eq.~(\ref{defcnc}),
\begin{eqnarray}
\label{t246}
{\mathcal T}_{\mathrm{CNC}}(246) 
&=& 2.37254~51620~38445~67035~68130~69148~85258~25756~18499~54254
\nonumber\\[1ex]
& & \;\;\; 97013~57806~20011~72404~62937~46020~32218~23862~67095~00004
\nonumber\\[1ex]
& & \;\;\; 69194~36541~28946~10390~15116~52595~90270~23975~58737~74256
\nonumber\\[1ex]
& & \;\;\; 23420~48480~95165~00802~19816~35378~76591~98589~60393~32102~8\,,
\end{eqnarray}
and
\begin{eqnarray}
\label{t247}
{\mathcal T}_{\mathrm{CNC}}(247) 
&=& 2.37254~51620~38445~67035~68130~69148~85258~25756~18499~54254
\nonumber\\[1ex]
& & \;\;\;  97013~57806~20011~72404~62937~46020~32218~23862~67095~00004
\nonumber\\[1ex]
& & \;\;\; 69194~36541~28946~10390~15116~52595~90270~23975~58737~74256
\nonumber\\[1ex]
& & \;\;\; 23420~48480~95165~00802~19816~35378~76591~98589~60393~32111~7\,.
\end{eqnarray}
The apparent convergence to roughly 200 decimals can be verified
against the right-hand side of~Eq.~(\ref{SimplifySumBailey}).
Of course, the right-hand side of~Eq.~(\ref{SimplifySumBailey}),
which involves only rationals, zeta functions and
the Euler--Mascheroni constant,
\[
\frac{17}{4} \, \zeta(4) - 4 \gamma\,\zeta(3) + \gamma^2 \, \zeta(2)\,,
\]
can be easily evaluated to 200 decimals
using known algorithms which are included in 
computer algebra systems (e.g.~\cite{Wo1988,disclaimer}).

The evaluation of the terms ${\bar a}(k)$ proceeds as follows.
For small index $k$, it is easy to write a recursion
relation relating ${\bar a}(k)$ and ${\bar a}(k+1)$ based on the 
(trivial) recursion for the $\psi$ function,
\begin{equation}
\psi(k+1) = \psi(k) + \frac{1}{k}\,.
\end{equation}
For large $k$, the asymptotic formula (\ref{PolyAsym})
can be used in order to calculate the $\psi$ function to
high precision. The point at which one may switch
from the recursion to the asymptotic method depends on how many
explicit values for Bernoulli numbers are available to the 
machine. We use the values for the first 60 Bernoulli numbers,
to 250 decimals, for our calculation.
We switch from one method to the other when the index $k$
of ${\bar a}(k)$ has reached a value of 500.

With 84~308 ${\bar a}(k)$ terms evaluated (out of which 
1364 by recursion and 82944 by the asymptotic method),
we evaluate the first 247 transforms with the results
presented above in Eqs.~(\ref{t246}) and (\ref{t247}).
As mentioned above,
if the terms of the series (\ref{SumBailey}) were
added on a term-by-term basis, then about $10^{205}$
would be required for an accuracy of 200 decimals in the 
final result. The reduction of this number to
roughly 84~000 corresponds to an acceleration of the convergence by
roughly 200 orders of magnitude.
We have also carried out, using enhanced precision, 
a calculation to 430 decimals,
involving about 500 CNC transformations and
arithmetic with 600 decimal figures. These evaluations
not only confirm the relation (\ref{SimplifySumBailey}) to high 
precision, but they also represent an accurate experimental
verification of the convergence properties of the 
delta transformation (\ref{defdelta}) in higher transformation order.
Moreover, it is observed that the rate of convergence of the
CNC transform results in a gain of approximately one significant
decimal figure per transformation, and that the rate of convergence 
remains constant over a wide range of transformation
orders (corresponding to linear convergence in the asymptotic region).
By contrast, the series (\ref{SimplifySumBailey}) is only logarithmically 
convergent. This corresponds to convergence
acceleration according to the definition~(\ref{ConvAccelDef}).

%
%
\section{Other Applications of the CNCT}
\label{sec_other}

We also mention the existing applications in the domain
of quantum electrodynamic bound-state 
calculations~(see~\cite{JeMoSo1999,Je1999,JeMoSo2001pra}. 
Another application concerns the quantum 
electrodynamic effective action (see~\cite{JeGiVaLaWe2001}).
The combined nonlinear-condensation transformation is
also applicable~\cite{Be1999master} to series of the form
\begin{equation}
\label{RpSer} 
R_p(x) \;=\; 
\sum_{k=0}^\infty \; \frac{x^{2k+1}}{(2k+1)^p} \,,
\end{equation}
\begin{equation}
\label{TpSer} 
T_p(x,b) \;=\; 
\sum_{k=0}^\infty \; \frac{1} {(2k+1)^p} \frac{ \cosh(2k+1)x}{\cosh(2k+1)b} \,,
\end{equation}
and
\begin{equation}
\label{UpSer} 
U_p(x,b) \;=\; 
\sum_{k=0}^\infty \; \frac{1}{(2k+1)^p}\frac{\cosh(2k+1)x}{\sinh(2k+1)b} \,. 
\end{equation}
Series of this type occur naturally in the context of plate
contact problems with periodic boundary conditions~\cite{DeKePaGl1984}.
The arguments $p$, $x$ and $b$ are real and positive for cases
of practical relevance. For $x \approx b$ and $p \approx 1$,
the series $T_p$ and $U_p$ are
very slowly convergent. In App.~A.2 of~\cite{Be1999master}
(pp.~105~ff.~{\em ibid.}), it
is demonstrated by way of numerical example that the CNCT is
able to efficiently accelerate the convergence of these series
in problematic parameter regions. 

In the numerical calculations, it is necessary to evaluate
terms with large index $k$. This can lead to numerical overflow 
because of the large arguments of the hyperbolic functions.
Clearly, representations such as
\begin{equation}
T_p(x,b) \;=\; 
\sum_{k=0}^\infty \; \frac{e^{(2k+1)(x-b)}}{(2k+1)^p} \;
\frac{1+e^{-2x(2k+1)}}{1+e^{-2b(2k+1)}}
\end{equation}
provide a solution for this problem.

Let us recall that considerable effort has been invested in the 
development of efficient numerical methods for the evaluation
of the series (\ref{RpSer}) -- 
(\ref{UpSer})~\cite{Ma1997,BaGa1995,BoDe1992,Ga1991,DeLiDe1990}.
These alternative methods make intensive use of special properties
of the series. They involve integral transformations and
infinite series over numerical integral~\cite{Ga1991}, and
they make use of
special properties of Legendre's chi-Function~\cite{BoDe1992}
which is related to the functions (\ref{RpSer}) --
(\ref{UpSer}).

We also
briefly mention that it is possible, in combining analytic results
obtained in~\cite{Pa1993,JePa1996,JeSoMo1997} 
with numerical techniques based on the
CNCT for the evaluation of complete hypergeometric functions, 
to evaluate the so-called Bethe logarithm in hydrogen to
essentially arbitrary precision. 
Specifically, we obtain -- for the
4P state -- the result
\begin{equation}
\ln k_0(4 {\mathrm P}) = -0.041~954~894~598~085~548~671~037(1)\,,
\end{equation}
which should be compared to other recent 
calculations~\cite{HaMo1985,DrSw1990,FoHi1993}.

%
%
\section{Conclusions}
\label{conclu}

We have discussed several applications of the
convergence acceleration methods introduced in Sec.~\ref{sec_cnct}:
in statistical physics (Sec.~\ref{sec_stat}),
in experimental mathematics (Sec.~\ref{sec_math}),
and other applications, mainly in the evaluation of
special functions (Sec.~\ref{sec_other}).
Specifically, it is observed that the 
combined nonlinear-condensation transformation
(CNCT) leads to an efficient calculational
scheme for the Lerch transcendent $\Phi$ given by 
Eq.~(\ref{lerch}). This special function provides a generalization
of several kinds of probability density functions which
are of significance for the statistical analysis e.g.
of DNA sequences (see Table~\ref{speccases}).
The comparatively fast and accurate evaluation of the 
statistical distributions with the CNCT is helpful with respect
to statistical procedures such as random number generation and statistical 
inference~\cite{AkSa2001}.

The high-precision calculation described in Sec.~\ref{sec_math} confirms that 
the rate of convergence of the CNC transforms is consistent 
with {\em linear} convergence (see Sec.~\ref{cnct_motivation})
in large transformation orders,
and in extended numerical precision (up to 430 digits).
In the absence of a rigorous proof~\cite{We1989} regarding the convergence
of the nonlinear sequence transformation (\ref{defdelta}),
which forms the second step of the CNCT (see Sec.~\ref{cnct_formulas}),
we attempt to gain numerical evidence for the (linear)
convergence of the CNCT in higher transformation orders,
and, in Sec.~\ref{sec_other}, for the wide applicability
of this algorithm in the context of the special functions
which are relevant to theoretical physics.

%
%
\section*{Acknowledgements}

The work by S.~V. Aksenov and M.~A. Savageau was supported in part by
U.S. Public Health Service Grant RO1-GM30054 from the National
Institutes of Health. U.D.J.~acknowledges
funding by the Deutsche Forschungsgemeinschaft (Nachwuchsgruppe
within the Sonderforschungsbereich 276).

\appendix
%
%
\section{Mathematica Program for Lerch's Transcendent}
\label{app}

The listing of the Mathematica program for the
calculation of Lerch's transcendent is shown in Fig.~\ref{mathprog}.   
The code displayed in Fig.~\ref{mathprog} is compact, but it lacks 
input argument checking. Also, there is no convenient control over the 
following parameters: (i) relative accuracy of the result, (ii) maximum 
number of allowed iterations, (iii) shift parameter $\beta$, and 
(iv) ``initial element'' parameter $n$. The code should be understood
as a template for more specialized implementations.
The program in Fig.~\ref{mathprog} 
takes advantage of Mathematica's implementation of the 
functional programming paradigm.

Necessarily, the Mathematica code shown in Fig.~\ref{mathprog}
is restricted to positive argument $z > 0$.
For the other arguments, the same restrictions exist as for the 
C program available from~\cite{AkJeURL}: 
$s,v$ must be real, the case of negative integer $v$
is excluded, and for negative non-integer $v$, $s$ is required 
to be integer.

\begin{figure}
\begin{verbatim}
1:  LerchPhiCNCT[z_,s_,v_]:=Module[{j,omega0,storeaj,num,den,o,res,factors},
2:  ajstep[inp_List]:=Module[{k,sum,bjk},
3:  k=inp[[1]];sum=inp[[2]];
4:  bjk=2^k z^(2^k (j+1)-1)/(v+2^k (j+1)-1)^s;sum+=bjk;k++;
5:  Return[{k,sum,bjk}]];
6:  recur[inp_List,pos_]:=Module[{loc},
7:  loc=inp;loc[[pos]]=loc[[pos+1]]-loc[[pos]] factors[[pos]];
8:  Return[loc]];
9:  sknstep[inp_List]:=Module[{i,omega,sn,skn,eps},
10: i=inp[[1]];omega=inp[[2]];sn=inp[[3]];skn=inp[[4]];eps=inp[[5]];
11: sn+=omega;AppendTo[storeaj,(-1)^i omega];
12: omega=(-1)^(i+1) If[EvenQ[i],0.5 (storeaj[[i/2+1]]-z^(i/2)/(v+i/2)^s),
    j=i+1;NestWhile[ajstep,{1,z^j/(v+j)^s,z^j/(v+j)^s},
    (#[[3]]/#[[2]]>10^(-2-acc))&][[2]]];
13: AppendTo[num,sn/omega];AppendTo[den,1/omega];
14: factors=Reverse[Which[i==0,{0},i==1,{0,1},True,
    Prepend[Table[(beta+n+i-1) (beta+n+i-2)/(beta+n+i+o-2)/(beta+n+i+o-3),
    {o,1,i}],0]]];
15: num=Fold[recur,num,Table[o,{o,Length[num]-1,1,-1}]];
16: den=Fold[recur,den,Table[o,{o,Length[den]-1,1,-1}]];
17: skn=RotateLeft[skn];eps=RotateLeft[eps];
18: skn[[2]]=num[[1]]/den[[1]];eps[[2]]=Abs[skn[[2]]-skn[[1]]];i++;
19: Return[{i,omega,sn,skn,eps}]];
20: acc=14;imax=100;beta=1;n=0;
21: j=0;omega0=NestWhile[ajstep,{1,z^j/(v+j)^s,z^j/(v+j)^s},
    (#[[3]]/#[[2]]>10^(-2-acc))&];num={};den={};storeaj={};
22: res=NestWhile[sknstep,{0,omega0[[2]],0,{0,0},{0,0}},
    Not[#[[1]]>imax-1||(#[[1]]>1&&(#[[5,2]]==0||(#[[5,2]]<#[[5,1]]&&
    Abs[2 (#[[5,1]])^2/(#[[5,1]]-#[[5,2]])/#[[4,2]]]<10^(-acc))))]&];
23: If[res[[1]]>imax-1,
    Print["Algorithm has not achieved relative accuracy of ",acc,
    " digits after maximum "res[[1]]," iterations."]];
24: Return[res[[4, 2]]]]
\end{verbatim}
\caption{A Mathematica program for the 
calculation of Lerch's transcendent using the CNC transformation.  
A line-by-line explanation is in the text.}
\label{mathprog}
\end{figure}

The overall structure 
of the program is the following: lines 2 -- 19 contain declarations of 
the three subroutines that are used iteratively during the calculation, 
and lines 20 -- 24 contain the main body. We explain the
algorithm underlying the program on a line-by-line basis.

\begin{itemize}
\item[1:] Declare the {\tt LerchPhiCNCT} function.
\item[2:] Declare the subroutine {\tt ajstep} that performs one 
iteration for calculation of the Van Wijngaarden transforms ${\bf A}_j$ 
according to Eqs.~(\ref{A2B}) and~(\ref{B2a}). The subroutine {\tt ajstep}
comprises the program lines 2 -- 5.
\item[3:] Read the current iteration counter $k$ and the 
temporary value of the partial sum of the series defining
${\bf A}_j$ [see Eq.~(\ref{A2B})] into local variables.
\item[4:] Calculate a new ${\bf b}_{k}^{(j)}$, update the temporary
value for the sum ${\bf A}_j$, and increment the counter $k$.
\item[5:] Return the result of one {\tt ajstep} iterate,
which is a list of the counter $k$, the current value of partial
sum defining ${\bf A}_j$, and the last calculated ${\bf b}_{k}^{(j)}$.
\item[6:] Declare the subroutine {\tt recur} that performs one iterative
recursive evaluation of the numerator and the denominator of the 
sequence transformation in Eq.~(\ref{SidTr}). The algorithm takes
advantage of the recursion relation in Eq.~(3.11) of~\cite{JeMoSoWe1999}.
The subroutine {\tt recur} comprises lines 6 -- 8 of the program.
\item[7:] Read the input into a local variable and apply the recurrence 
relation.
\item[8:] Return the result of {\tt recur}, which is an 
updated list {\tt loc}.
\item[9:] Declare the subroutine {\tt sknstep} which performs  
iteratively the calculation of the CNC transforms. 
The subroutine {\tt recur} comprises lines 9 -- 19.
\item[10:] Read the following data into local
variables: (i) the current iteration counter $i$, 
(ii) the Van Wijngaarden term (also remainder estimate) 
$\omega_{i-1} = (-1)^i \, {\bf A}_i$, 
(iii) the partial sum of the Van Wijngaarden series $s_i$ 
(which is the input data for the delta transforms), 
(iv) the last two CNC transforms calculated
[$\delta^{(0)}_{i-2} \equiv {\cal S}_{i-2}^{(0)}$
and $\delta^{(0)}_{i-1}$], 
(v) the difference between CNC transforms of 
successive iterations $\epsilon_{i-2}$ and
$\epsilon_{i-1}$ where $\epsilon_{k} \equiv |\delta^{(0)}_{k}  -
\delta^{(0)}_{k-1}|$. The definitions of
these symbols can be found in Eqs.~(\ref{VWSerTran}) -- (\ref{B2a}),
(\ref{SidTr}), (\ref{defdelta}), and (\ref{repl}). 
\item[11:] Increment the partial sum $s_i$ with the 
next term $(-1)^i \, {\bf A}_i$ of the Van Wijngaarden series,
and store the current ${\bf A}_i$ in the list of variables {\tt storeaj}
for later use in the recurrence relation Eq.~(\ref{compute}).
\item[12:] Calculate the next term 
$\omega_i = (-1)^{i+1} \, {\bf A}_{i+1}$ of the 
Van Wijngaarden series: for
even index $i$, use Eq.~(\ref{compute}).
For odd $i$, calculate via a call to the subroutine {\tt ajstep}. 
The last argument in the ``nested call'' to the Mathematica
built-in function {\tt NestWhile} is a termination criterion: the 
${\bf b}_{k}^{(j)}$ are added on a term-by-term basis until
the ratio of the third element of the 
list of the results of {\tt ajstep} (which is
the value of the current ${\bf b}_{k}^{(j)}$) 
and the second element of that list 
(which is the value of the current partial sum of
${\bf A}_j$) is less than the specified accuracy.  
This termination criterion is justified for the case of the Lerch 
transcendent by Eq.~(\ref{criterion}).
\item[13:] Calculate and append new values in the lists
that store the numerator and denominator of the 
$\delta^{(0)}_i \equiv {\cal S}_{i}^{(0)}$.
The length of these lists corresponds to the current order 
of the CNC transformation $i$.
\item[14:] Calculate a list of prefactors used for the 
recursive evaluation of the denominator and numerator 
of the delta transformations
$\delta^{(0)}_i$ according to Eq.~(3.11) of~\cite{JeMoSoWe1999}.
\item[15:] Calculate the numerator of the delta transforms
recursively and store the results in the list {\tt num}.
\item[16:] Calculate the denominator recursively
and store results in the list {\tt den}.
\item[17:] Swap the elements of the 
lists {\tt skn} and {\tt eps} which contain
the last two delta transforms and the last two 
absolute differences; the purpose is to let the ``old'' values 
be the first elements, in order to store the ``newly calculated''
values as the second elements in a later step.
\item[18:] Calculate a ``new'' delta transform and a new
absolute difference {\tt eps} and store 
them as second elements of {\tt skn} and {\tt eps}.
\item[19:] Return the result of {\tt sknstep}.
It consists of a ``new'' iteration counter {\tt i}, 
a ``new'' Van Wijngaarden term {\tt omega}, 
a ``new'' partial sum of the Van Wijngaarden series {\tt sn}, 
a list {\tt skn} of the last two CNC transforms, 
and a list {\tt eps} of the last two differences between CNC transforms 
of successive iterations.
\item[20:] Start the main program.
Set the default values for the desired relative accuracy of 
the result, for the maximum allowed order of the delta 
transformation, for the shift parameter 
$\beta = 1$, and for the ``initial element'' parameter $n=0$ 
which enters the CNC transform 
$\delta_{k}^{(n)} \equiv {\cal S}_{k}^{(n)}$ in Eq.~(\ref{SidTr}).
\item[21:] Calculate the initial remainder 
estimate $\omega$ via a nested call 
to {\tt ajstep} (in analogy to line 12) and initialize other variables.
\item[22:] Iteratively compute the CNC transforms
using {\tt sknstep}. Again, in analogy to line 12,
the second argument of the {\tt NestWhile} function is the list of initial 
values, and the third argument is the termination criterion. 
The evaluation of the CNC transforms stops if either the 
maximum number of iterations is 
exceeded (first element of the output of {\tt sknstep}), or if 
the current difference between the transforms is zero (fifth element), 
or if an error estimate for the last evaluated delta transform 
is less than the specified accuracy. This error estimate is inspired
by the fact that the delta transforms are observed to converge
linearly in good approximation [see Eq.~(\ref{ConvAccelDef2}) and
Sec.~\ref{cnct_motivation} as well as Sec.~\ref{sec_math}].
The remainder estimate is 
evaluated using the $\epsilon_{i-1}$ and 
$\epsilon_{i}$ (calculated from the fifth argument 
to {\tt sknstep}) and the current delta transform 
(which is second element of the fourth argument 
to {\tt sknstep}).
\item[23:] Print an error message if the maximum number of 
iterations is exceeded.
\item[24:] Return the result of {\tt LerchPhiCNCT}, which is 
the second element of the
fourth argument of {\tt sknstep} and corresponds
to the last delta transform evaluated.
\end{itemize}
To conclude this Appendix, we provide a more detailed
explanation for the error estimate used in program line 22.

The output of the CNC transformation is a sequence
of approximants ${\cal T}_n \equiv {\cal T}_{\rm CNC}(n)$ 
[see Eq.~(\ref{defcnc})] that converge
to the value of $\Phi(z,s,v)$. We define the ratio of two
consecutive differences of approximants as
\begin{equation}
\label{xx}
x_n = \left| \frac{{\cal T}_{n} - {\cal T}_{n-1}}
  {{\cal T}_{n-1} - {\cal T}_{n-2}} \right|\,.
\end{equation}
As discussed in Sec.~\ref{sec_math}, the CNC transforms can be expected 
to converge geometrically in higher transformation orders.  
Let $\bar \rho$ be defined as [see also Eq.~(\ref{LimRatSeq})]
\begin{equation}
\bar \rho = \lim_{n \to \infty} 
\frac {{\cal T}_{n} - {\cal T}_\infty} {{\cal T}_{n-1} - {\cal T}_\infty}\,,
\end{equation}
with ${\cal T}_\infty = \Phi(z,s,v)$.
Clearly, the $x_n$ will provide a good
estimate for $\bar \rho$ at large $n$.
A good estimate for the truncation error
${\cal T}_{n} - {\cal T}_\infty$ can thus be obtained
by summing the geometric series
$\sum_{k=1}^\infty {\bar \rho}^k \, |{\cal T}_{n} - {\cal T}_{n-1}|$ 
where our best estimate for ${\bar \rho}$ is ${\bar \rho} \approx x_n$.
We may therefore use the following convergence criterion to terminate
the calculation of the CNC transforms:
\begin{equation}
\label{convcrit}
\frac{2}{x_n} \, \left[ \frac{1}{1-x_n} \,
\left| \frac{{\cal T}_{n} - {\cal T}_{n-1}}{{\cal T}_{n}} \right| \right] < 
{\tt acc}\,.
\end{equation}
Here, ${\tt acc}$ is the specified desired {\em relative}
accuracy of the result. The factor $2/x_n$ in~(\ref{convcrit})
is a heuristic ``safeguard factor'' introduced with the
notion of avoiding a premature termination of the calculation
of successive transforms in the problematic case of
two consecutive transforms accidentally assuming values
very close to each other. Such a situation may arise {\em before}
the asymptotic, geometric convergence sets in.
The term in square brackets in (\ref{convcrit}) represents
the remainder estimate based on the geometric model $x_n \approx 
{\bar \rho}$ [see Eq.~(\ref{xx})]. Rewriting (\ref{convcrit}),
we obtain the termination criterion used in line 22 of the
program in Fig.~\ref{mathprog}.

Finally, we would like to point out that the above 
implementation relies on only one algorithm (CNCT), 
and that one cannot expect to
obtain optimal performance in all parameter regions, let alone
analytic continuations for those cases where the power series (\ref{lerch})
diverges. As regards the evaluation
of special functions with very large (excessive)
parameter values, it is known that asymptotic expansions can
provide optimal methods of evaluation (see, e.g.~\cite{Mo1974b}). 
These are not implemented in the code in Fig.~\ref{mathprog}.
For negative $z$, in our case
some improvement can be obtained by accelerating
the rate of convergence of the resulting alternating series
by the direct application~\cite{We1989,AkJeURL} of delta 
transforms to the power series (\ref{lerch});
this avoids the potentially time-consuming Van Wijngaarden step
(\ref{VWSerTran}) -- (\ref{B2a}) of the CNCT. However,
the more problematic case of positive $z$, which results in
a {\em non}alternating series, requires the (full)
combined nonlinear-condensation transformation implemented by
the algorithm in Fig.~\ref{mathprog}.


\begin{thebibliography}{10}
\bibitem{disclaimer}
Certain commercial equipment, instruments, or materials are identified in this
  paper to foster understanding. Such identification does not imply
  recommendation or endorsement by the National Institute of Standards and
  Technology, nor does it imply that the materials or equipment identified are
  necessarily the best available for the purpose.
\end{thebibliography}
\end{document}